\documentclass[12pt]{article}
\usepackage{a4,amsmath,amsthm,amssymb,enumerate,eucal,mathrsfs,graphicx,
  subfig,xcolor,array}
\title{\textbf{Quasigraphs and skeletal partitions}}%
\author{Tom\'{a}\v{s} Kaiser$^{\:1}$\\
  Petr Vr\'{a}na$^{\:1}$}
\date{}

\newtheorem{theorem}{Theorem}
\newtheorem{lemma}[theorem]{Lemma}

\newtheorem{claim}{Claim}

\newtheorem{corollary}[theorem]{Corollary}

\newtheorem{observation}[theorem]{Observation}

\newcommand\size[1] {\left|{#1}\right|}
\newcommand\Set[2] {\left\{{#1}:\,{#2}\right\}}
\newcommand\Setx[1] {\left\{{#1}\right\}}

\newcommand{\MM}{\mathscr M}%
\newcommand{\PP}{\mathcal P}%
\newcommand{\QQ}{\mathcal Q}%
\newcommand{\RR}{\mathcal R}%
\newcommand{\SSS}{\mathcal S}%
\newcommand{\TT}{\mathcal T}%
\newcommand{\lseq}[1]{\tilde{\mathbb{P}}^{#1}}
\newcommand{\ptn}[3]{\PP^{#1}_{#2,#3}}

\newcommand{\sw}[2]{{#1}^{({#2})}}

\newcommand{\worse}{\sqsubset}
\newcommand{\worseeq}{\sqsubseteq}
\newcommand{\worsex}[2]{\sqsubset_{(#1,#2)}}
\newcommand{\worseeqx}[2]{\sqsubseteq_{(#1,#2)}}
\newcommand{\eq}{\equiv}
\newcommand{\eqx}[2]{\eq_{(#1,#2)}}
\newcommand{\pairs}{\TT}
\newcommand{\stopx}{\text{\textsc{stop}}}
\newcommand{\termx}{\text{\textsc{terminate}}}
\newcommand{\sig}[1]{\mathbb S^{#1}}

\newcommand{\fig}[1]{\includegraphics[page=#1]{skeletal-figs.pdf}}%
\newcommand{\sfig}[2]{\subfloat[#2]{\fig{#1}}}%

\newcommand{\sfigtop}[2]{\newbox{\pic}\sbox{\pic}{\fig{#1}}%
  \subfloat[#2]{\vbox to\ht\base{\hbox to\wd\pic{\usebox\pic}}}}
\newcommand{\hf}{\hspace*{0pt}\hspace*{\fill}\hspace*{0pt}}

\newcommand{\claimproofend}{\hspace*{.1mm}\hspace{\fill}}
\newenvironment{claimproof}{}{\claimproofend\par\vspace{2mm}}
\newtheorem{xcasehdr}{Case}
\newtheorem{xsubcasehdr}{Subcase}[xcasehdr]
\newtheorem{xsubsubcasehdr}{Subsubcase}[xsubcasehdr]

\newenvironment{xcase}[1]%
{\vspace{-1mm}\par\noindent\xcasehdr{#1}\upshape\vspace{2mm}\par\noindent}%
{}

\newenvironment{xsubcase}[1]%
{\vspace{-1mm}\par\noindent\xsubcasehdr{#1}\upshape
  \vspace{2mm}\par\noindent}%
{\hspace*{0mm}\par}

{\vspace{-1mm}\par\noindent\xsubsubcasehdr{#1}\upshape
  \vspace{2mm}\par\noindent}%
{\hspace*{0mm}\par}

\begin{document}
\maketitle
\footnotetext[1]{Department of Mathematics and European Centre of
  Excellence NTIS (New Technologies for the Information Society),
  University of West Bohemia, Univerzitn\'{\i}~8, 306~14~Plze\v{n},
  Czech Republic. Supported by projects 17-04611S and 20-09525S of the
  Czech Science Foundation. E-mail:
  \{\texttt{kaisert,vranap}\}\texttt{@kma.zcu.cz}.}

\begin{abstract}
  We give a new proof of the Skeletal Lemma, which is the main
  technical tool in our paper on Hamilton cycles in line
  graphs~[T. Kaiser and P.~Vr\'{a}na, Hamilton cycles in 5-connected
  line graphs, European J. Combin. 33 (2012), 924--947]. It
  generalises results on disjoint spanning trees in graphs to the
  context of 3-hypergraphs. The lemma is proved in a slightly stronger
  version that is more suitable for applications. The proof is
  simplified and formulated in a more accessible way.
\end{abstract}


\section{Introduction}
\label{sec:introduction}

The main tool used in our work on Hamilton cycles in line
graphs~\cite{KV-hamilton} is a result called `Skeletal
Lemma'~\cite[Lemma 17]{KV-hamilton}. It deals with quasigraphs in
3-hypergraphs (see below for definitions) and is related to Tutte's
and Nash-Williams' characterisation of graphs with two disjoint
spanning trees.

In our recent paper~\cite{KV-ess-9}, we need to use the lemma in a
slightly stronger form that unfortunately does not follow from the
formulation given in the paper. Instead of pointing out the necessary
modifications to the long and complicated proof, we decided to use
this opportunity to rewrite the proof completely, trying to formulate
it in a way as conceptually simple as we can. That is the purpose of
the present paper which is a companion paper to~\cite{KV-ess-9}. In
addition, the present paper aims to give the full proof in detail,
even in parts where the argument in \cite{KV-hamilton} is somewhat
sketchy.

The structure of the paper is as follows. In Section~\ref{sec:quasi},
we review the basic notions related to quasigraphs, the structures
forming a central concept of our proof. In
Section~\ref{sec:anti-conn}, we develop the basic properties of the
notion of connectivity and especially anticonnectivity of a quasigraph
on a set of vertices. This allows us to define, for any quasigraph, a
sequence of successively more and more refined partitions of the
vertex set that serves as a measure of `quality' of the
quasigraph. This is done in
Section~\ref{sec:sequence}. Section~\ref{sec:skeletal} gives the proof
of the main result, a stronger version of the Skeletal Lemma
(Theorem~\ref{t:enhancing}). Finally, in Section~\ref{sec:bad}, we
infer the result we need for the above mentioned application
in~\cite{KV-ess-9} (Theorem~\ref{t:no-bad}).


\section{Quasigraphs}
\label{sec:quasi}

A \emph{$3$-hypergraph} is a hypergraph whose hyperedges have size $2$
or $3$. Throughout this paper, let $H$ be a 3-hypergraph. A
\emph{quasigraph} in $H$ is a mapping $\pi$ that assigns to each
hyperedge $e$ of $H$ either a subset of $e$ of size 2, or the empty
set. The hyperedges $e$ with $\pi(e)\neq \emptyset$ are said to be
\emph{used} by $\pi$. (See Figure~\ref{fig:quasi} for an
illustration.) 
A quasigraph $\pi$ in $H$ restricts to a quasigraph in
any subhypergraph $H'$ of $H$; to avoid excessive notation, we will
use $\pi$ to denote the corresponding quasigraph in $H'$ as well.

\begin{figure}
  \begin{center}
    \hf\subfloat[]{\fig3}\hf\subfloat[]{\fig4}\hf
  \end{center}
  \caption{(a) A 3-hypergraph $H$. Hyperedges of size $3$ are
    depicted as three lines meeting at a point without a vertex
    mark. (b) A quasigraph $\pi$ in $H$. For each hyperedge $e$ used
    by $\pi$, the pair $\pi(e)$ is shown using one or two bold
    lines, depending on the size of $e$. The other lines are shown
    dotted for greater contrast.}
  \label{fig:quasi}
\end{figure}

Given a quasigraph $\pi$ in $H$, we let $\pi^*$ denote the graph on
$V(H)$, obtained by considering the pairs $\pi(e')$ as edges whenever
$\pi(e')\neq\emptyset$ ($e'\in E(H)$). If $\pi^*$ is a forest, then
$\pi$ is \emph{acyclic}. If $\pi^*$ is the union of a cycle and a set
of isolated vertices, then $\pi$ is a \emph{quasicycle}. A
3-hypergraph $H$ is \emph{acyclic} if there exists no quasicycle in
$H$.

If $e$ is a hyperedge of $H$, then we define $\pi-e$ as the quasigraph
which satisfies $(\pi-e)(e) = \emptyset$, and coincides with $\pi$ on
all hyperedges other than $e$. If $e$ is a hyperedge not used by
$\pi$, and if $u,v\in e$, then $\pi+(uv)_e$ is the quasigraph that
coincides with $\pi$ except that its value on $e$ is $uv$ rather than
$\emptyset$.

The \emph{complement} $\overline\pi$ of $\pi$ is the subhypergraph of $H$
(on the same vertex set) consisting of the hyperedges not used by
$\pi$.

Let $\PP$ be a partition of a set $X\subseteq V(H)$. We say that $\PP$
is \emph{nontrivial} if $\PP\neq\Setx{X}$. For $Y\subseteq X$, the
partition $\PP[Y]$ of $Y$ \emph{induced} by $\PP$ has all nonempty
intersections $P\cap Y$, where $P\in\PP$, as its classes. If a
hyperedge $e$ of $H$ intersects at least two classes of $\PP$, then it
is said to be \emph{$\PP$-crossing}.

Assume now that $X = V(H)$, i.e., $\PP$ is a partition of $V(H)$.  If
$e\in E(H)$, then $e/\PP$ is defined as the set of all classes of
$\PP$ intersected by $e$. The hypergraph $H/\PP$ has vertex set $\PP$
and its hyperedges are all the sets of the form $e/\PP$, where $e$ is
a $\PP$-crossing hyperedge of $H$. Thus, $H/\PP$ is a 3-hypergraph. A
quasigraph $\pi/\PP$ in this hypergraph is defined by setting, for
every $\PP$-crossing hyperedge $e$ of $H$,
\begin{equation*}
  (\pi/\PP)(e/\PP) =
  \begin{cases}
    \pi(e)/\PP & \text{if $\pi(e)$ is $\PP$-crossing,}\\
    \emptyset & \text{otherwise.}
  \end{cases}
\end{equation*}

We extend the above notation and write, e.g., $uv/\PP$ for the set of
classes of $\PP$ intersecting $\Setx{u,v}$, where $u,v\in V(H)$.

By the above definitions, the complement of $\pi/\PP$ is the
subhypergraph of $H/\PP$ consisting of all the hyperedges $e/\PP$ such
that $\pi(e)$ is contained in some class of $\PP$ (including
$\pi(e) = \emptyset$). We often consider quasigraphs $\gamma$ in
$\overline{\pi/\PP}$ (typically, such a $\gamma$ is a quasicycle). An
example is given in Figure~\ref{fig:complement}.

\begin{figure}
  \begin{center}
    \hf\subfloat[]{\fig5}\hf\subfloat[]{\fig6}\hf
    \caption{(a) A quasigraph $\pi$ in a $3$-hypergraph $H$ and a
      partition $\PP$ of $V(H)$. The classes of the partition are
      shown in grey. (b) A quasicycle $\gamma$ in the complement of
      $\pi/\PP$ in the hypergraph $H/\PP$. Note that the vertex set of
      this hypergraph is $\PP$.}
    \label{fig:complement}
  \end{center}
\end{figure}


\section{(Anti)connectivity}
\label{sec:anti-conn}

In this section, we define and explore the notions of components and
anticomponents of a quasigraph (on a set of vertices) that are
completely essential for our arguments. We refer to
Figure~\ref{fig:anti} for an illustration of these notions.

Recall that $H$ denotes a 3-hypergraph. Let $\pi$ be a quasigraph in
$H$ and $X\subseteq V(H)$. We say that $\pi$ is \emph{connected on
  $X$} if the induced subgraph of $\pi^*$ on $X$ is connected. The
\emph{components} of $\pi$ on $X$ are defined as the vertex sets of
the connected components of the induced subgraph of $\pi^*$ on $X$.

We say that $\pi$ is \emph{anticonnected on $X$} if for each
nontrivial partition $\RR$ of $X$, there is an $\RR$-crossing
hyperedge $f$ of $H$ such that $\pi(f)$ is a subset of one of the
classes of $\RR$ (possibly $\pi(f)=\emptyset$). If we need to refer to
the hypergraph $H$, we say that $\pi$ is \emph{anticonnected on $X$ in
  $H$}.

Observe that $\pi$ is both connected and anticonnected on any set
consisting of a single vertex.

\begin{figure}
  \begin{center}
    \fig7
  \end{center}
  \caption{A quasigraph $\pi$ in $H$ and a set $X \subseteq V(H)$
    (shown grey). The quasigraph $\pi$ is anticonnected on $X$ and has
    four components on $X$.}
  \label{fig:anti}
\end{figure}

\begin{lemma}\label{l:union}
  Let $\pi$ be a quasigraph in (a $3$-hypergraph) $H$ and $X,Y$
  subsets of $V(H)$ such that $\pi$ is anticonnected on $X$ and
  $Y$. Then $\pi$ is anticonnected on $X\cup Y$ whenever one of the
  following holds:
  \begin{enumerate}[\quad(i)]
  \item $X$ and $Y$ intersect, or
  \item there is a hyperedge $h$ of $H$ intersecting both $X$ and $Y$,
    such that $\pi(h)$ is a subset of $X$ or $Y$ (possibly
    $\pi(h)=\emptyset$).
  \end{enumerate}
\end{lemma}
\begin{proof}
  Let $\RR$ be a nontrivial partition of $X\cup Y$. We find for $\RR$
  the hyperedge whose existence is required by the definition of
  anticonnectedness.

  Suppose first that $\RR[X]$ is nontrivial. Since $\pi$ is
  anticonnected on $X$, there is a hyperedge $f$ of $H$ such that $f$
  intersects at least two classes of $\RR[X]$ and one of them contains
  $\pi(f)$. Thus, $f$ intersects at least two classes of $\RR$ and one
  of them contains $\pi(f)$.

  We can thus assume, by symmetry, that both $\RR[X]$ and $\RR[Y]$ are
  trivial. This implies that $\RR=\Setx{X,Y}$, so $X$ and $Y$ are
  disjoint. In this case, the hyperedge $h$ from (ii) has the required
  property.
\end{proof}

By Lemma~\ref{l:union}, the maximal sets $Y\subseteq X$ such that
$\pi$ is anticonnected on $Y$ partition $X$. We call them the
\emph{anticomponents} of $\pi$ on $X$.

\begin{lemma}\label{l:no-change}
  Let $\pi$ and $\rho$ be quasigraphs in $H$ and $Y$ be a subset of
  $V(H)$ such that $\pi(e) = \rho(e)$ for every hyperedge $e$ of $H$
  with $\size{e\cap Y} \geq 2$. Then $\pi$ is anticonnected on $Y$ if
  and only if $\rho$ is anticonnected on $Y$.
\end{lemma}
\begin{proof}
  Suppose that $\pi$ is anticonnected on $Y$ and let $\RR$ be a
  nontrivial partition of $Y$. Consider a hyperedge $f$ of $H$ such
  that $f$ intersects two classes of $\RR$ and one of them contains
  $\pi(f)$. By the assumption, $\rho(f) = \pi(f)$, so the same holds
  for $\rho$ in place of $\pi$. Since $\RR$ is arbitrary, $\rho$ is
  anticonnected on $Y$. The lemma follows by symmetry.
\end{proof}

We prove several further lemmas that describe some of the basic
properties of (anti)con\-nect\-i\-vi\-ty of quasigraphs.

\begin{lemma}\label{l:H-e}
  Let $\pi$ be a quasigraph in $H$, $X\subseteq V(H)$ and $e$ a
  hyperedge of $H$ with $\size{e\cap X}\leq 1$. If $\pi$ is
  anticonnected on $X$ in $H$, then $\pi$ is anticonnected on $X$ in
  $H-e$.
\end{lemma}
\begin{proof}
  Let $\RR$ be a nontrivial partition of $X$. Since $\pi$ is
  anticonnected on $X$ in $H$, there is a hyperedge $f$ of $H$
  intersecting at least two classes of $\RR$, one of which contains
  $\pi(f)$. The hyperedge $f$ is distinct from $e$ as
  $\size{e\cap X} \leq 1$. Thus, $f \in E(H-e)$. Since $\RR$ is
  arbitrary, $\pi$ is anticonnected on $X$ in $H-e$.
\end{proof}

\begin{lemma}\label{l:sub-add}
  Let $\pi$ be a quasigraph in $H$ and $Y\subseteq X$ subsets of
  $V(H)$. Suppose that $e$ is a hyperedge of $H$ not used by $\pi$ and
  containing vertices $u,v\in Y$. Define $\rho$ as the quasigraph $\pi
  + (uv)_e$. The following holds:
  \begin{enumerate}[\quad(i)]
  \item if $\pi$ is anticonnected on $X$ and $\rho$ is anticonnected
    on $Y$, then $\rho$ is anticonnected on $X$,
  \item if $\pi$ is connected on $X$, then so is $\rho$.
  \end{enumerate}
\end{lemma}
\begin{proof}
  We prove (i). Consider an arbitrary partition $\RR$ of $X$. We aim
  to show that there is a hyperedge $f$ of $H$ such that $f$
  intersects two classes of $\RR$ and $\rho(f)$ is contained in one of
  them. This is certainly true if $\RR[Y]$ is nontrivial, since $\rho$
  is assumed to be anticonnected on $Y$. Thus, we may assume that $Y$
  is contained in a class of $\RR$.

  Since $\pi$ is anticonnected on $X$, there is a hyperedge $h$ of $H$
  such that $h$ intersects two classes of $\RR$ and $\pi(h)$ is
  contained in one of them. We set $f := h$. If $h\neq e$, this choice
  works because $\rho(h) = \pi(h)$. If $h=e$, then $\rho(h)$ is
  contained in $Y$ and therefore in a class of $\RR$. This concludes
  the proof of (i).

  Part (ii) follows directly from the fact that $\pi^*$ is a subgraph
  of $\rho^*$, and therefore the induced subgraph of $\pi^*$ on $X$ is
  a subgraph of the induced subgraph of $\rho^*$ on $X$.
\end{proof}

\begin{lemma}\label{l:sub-remove}
  Let $\pi$ be a quasigraph in $H$ and $Y\subseteq X$ subsets of
  $V(H)$. Suppose that $e$ is a hyperedge of $H$ with $\pi(e)\subseteq
  Y$. Define $\sigma$ as the quasigraph $\pi-e$. It holds that
  \begin{enumerate}[\quad(i)]
  \item if $\pi$ is anticonnected on $X$, then so is $\sigma$,
  \item if $\pi$ is connected on $X$ and $\sigma$ is connected on $Y$,
    then $\sigma$ is connected on $X$.
  \end{enumerate}
\end{lemma}
\begin{proof}
  We prove (i). Suppose, for contradiction, that $\sigma$ is not
  anticonnected on $X$.  By the definition, there is a partition
  $\SSS$ of $X$ such that for all hyperedges $f$ of $H$ intersecting
  at least two classes of $\SSS$, $\sigma(f)$ intersects two classes
  of $\SSS$ as well. On the other hand, since $\pi$ is anticonnected
  on $X$ and $\pi$ has the same value as $\sigma$ on every hyperedge
  other than $e$, it must be that the hyperedge $e$ intersects two
  classes of $\SSS$ (and $\pi(e)$ is contained in one class). Since
  $\sigma(e)=\emptyset$, we obtain a contradiction.

  Next, we prove (ii). Note that $\sigma^*$ is a subgraph of
  $\pi^*$. Since $\sigma$ is connected on $Y$, so is $\pi$.  We show
  that $\sigma$ is connected on $X$. Let $\pi^*_X$ be the induced
  subgraph of $\pi^*$ on $X$, and let $\pi(e) = \Setx{u,v}$. We need
  to prove that any two vertices in $X$ are joined by a walk in the
  induced subgraph of $\sigma^*$ on $X$, which equals
  $\pi^*_X-uv$. This is easy from the fact that $\pi^*_X$ is
  connected, and that the edge $uv$ may be replaced in any walk by a
  path from $u$ to $v$ in the induced subgraph of $\sigma^*$ on $Y$
  (which is connected).
\end{proof}

Let us now define two notions that will play a role when we introduce
the sequence of a quasigraph in Section~\ref{sec:sequence}. (See
Figure~\ref{fig:bridge} for an illustration.) Suppose that
$X\subseteq V(H)$ such that the quasigraph $\pi$ is both connected and
anticonnected on $X$. Let $e$ be a hyperedge with
$\size{e\cap X} = 2$.

We say that $e$ is an \emph{$X$-bridge} (with respect to $\pi$) if $e$
is used by $\pi$, $\pi(e) \subseteq X$, and $\pi-e$ is not connected
on $X$ in $H-e$. Similarly, $e$ is an \emph{$X$-antibridge} (with
respect to $\pi$) if $e$ is not used by $\pi$ and $\pi$ is not
anticonnected on $X$ in $H-e$.

\begin{figure}
  \begin{center}
    \fig8
  \end{center}
  \caption{A quasigraph $\pi$ in a $3$-hypergraph $H$ and a set
    $X\subseteq V(H)$ (shown grey) such that $\pi$ is both connected
    and anticonnected on $X$. The hyperedge $e$ is an $X$-bridge with
    respect to $\pi$, while $f$ is an $X$-antibridge with respect to
    $\pi$.}
  \label{fig:bridge}
\end{figure}


\section{The plane sequence of a quasigraph}
\label{sec:sequence}

Let $H$ be a $3$-hypergraph. Throughout our arguments, we will work
with partitions of $V(H)$. Given two such partitions $\PP$ and $\QQ$,
we say that $\PP$ \emph{refines} $\QQ$ (and write $\PP \leq \QQ$) if
each class of $\PP$ is a subset of some class of $\QQ$. If $\PP$
refines $\QQ$ and $\PP\neq\QQ$, we write $\PP < \QQ$.

Let $\pi$ be a quasigraph in $H$. In~\cite{KV-hamilton}, we associate
with $\pi$ a sequence of partitions of $V(H)$. In the present paper,
we proceed similarly, but for technical reasons, we need to extend the
original definition to involve a two-dimensional analogue of a
sequence. A \emph{plane sequence} is a family
$(\PP_{i,j})_{i,j\geq 0}$ of partitions of $V(H)$.

It will be convenient to consider the lexicographic order $\leq$ on
pairs of nonnegative integers: $(i,j) \leq (i',j')$ if either
$i < i'$, or $i = i'$ and $j \leq j'$. This is extended in the natural
way to the set
\begin{equation*}
  \pairs = \Set{(i,j)}{0\leq i < \infty, 0\leq j\leq\infty}\cup\Setx{(\infty,\infty)}.
\end{equation*}
For instance, $(1,\infty) < (2,0) < (\infty,\infty)$. This is a
well-ordering on the set $\pairs$, which allows us to perform
transfinite induction over $\pairs$ (cf.~\cite{Me}).

The \emph{(plane) sequence of $\pi$}, denoted by $\lseq\pi$, consists of
partitions $\ptn\pi i j$ of $V(H)$, where $(i,j)\in\pairs$. We let
$\ptn\pi00$ be the trivial partition $\Setx{V(H)}$. If $j \geq 1$ and
$\ptn\pi i {j-1}$ is defined, then we let
\begin{equation*}
  \ptn\pi i j =
  \begin{cases}
    \Set{K}{K\text{ is a component of $\pi$ on some
        $X\in\ptn\pi i {j-1}$}} & \text{if $j$ is
      odd,}\\
    \Set{K}{K\text{ is an anticomponent of $\pi$ on some
        $X\in\ptn\pi i {j-1}$}} & \text{if $j$ is even.}
  \end{cases}
\end{equation*}
See Figure~\ref{fig:seq} for an example.

\begin{figure}
  \centering\sfig{24}{}\hf\sfig{25}{}\\
  \sfig{26}{}\hf\sfig{23}{}
  \caption{(a) A quasigraph $\pi$ (bold) in a $3$-hypergraph $H$. (b)
    The single class of $\ptn\pi00$ (lighter gray) and the classes of
    $\ptn\pi01$ (darker gray). (c) The classes of $\ptn\pi01$
    (lighter) and $\ptn\pi02$ (darker). (d) All the partitions
    $\ptn\pi00,\dots,\ptn\pi03$ (lightest to darkest gray). In this
    case, $\ptn\pi03 = \ptn\pi0\infty$.}
  \label{fig:seq}
\end{figure}


So far, this yields the partitions $\ptn\pi00,\ptn\pi01,\dots$. We
first notice that since $H$ is finite, there is some $j_0$ such that
\begin{equation*}
  \ptn\pi 0 {j_0} = \ptn\pi 0 {j_0 + 1} = \dots,
\end{equation*}
and we set $\ptn\pi 0 \infty$ equal to $\ptn\pi 0 {j_0}$. We will use
an analogous definition to construct $\ptn\pi i \infty$ for $i > 0$
when $\ptn\pi i 0,\ptn\pi i 1,\dots$ will have been defined. (See
Figure~\ref{fig:layout} for a schematic illustration.)

\begin{figure}
  \centering\fig{10}
  \caption{The order of partitions in the construction of a plane
    sequence of a quasigraph.}
  \label{fig:layout}
\end{figure}

By the construction, $\ptn\pi 0 \infty$ has the property that $\pi$ is
both connected and anticonnected on each of its classes. We call any
such partition of $V(H)$ \emph{$\pi$-solid}.

The definition of the plane sequence of $\pi$ will be completed once
we define $\ptn\pi i 0$ for all $i\geq 1$. Thus, let $i\geq 1$ be
fixed, and suppose that the partition $\PP:=\ptn\pi {i-1} \infty$ is
already defined.

Let $A,B\in\PP$. The \emph{exposure step} of the pair $AB$ is the
least $(s,t)$ (with respect to the ordering defined above) such that
$A$ and $B$ are contained (as subsets) in different classes of
$\ptn\pi s t$. Similarly, for a pair of vertices $u,v$ of $H$
contained in different classes of $\PP$, the exposure step of the pair
$uv$ is the least $(s,t)$ such that $u$ and $v$ are contained in
different classes of $\ptn\pi s t$.

Suppose that $\gamma$ is a quasicycle in $H/\PP$. The \emph{exposure
  step} of $\gamma$ is the least exposure step of $\gamma(e/\PP)$,
where $e$ ranges over all hyperedges of $H$ such that $e/\PP$ is used
by $\gamma$. If the exposure step of $\gamma$ is $(s,t)$, we also say
that $\gamma$ is \emph{exposed} at (step) $(s,t)$. We say that a
hyperedge $e$ of $H$ is a \emph{leading hyperedge} of $\gamma$ if
$e/\PP$ is used by $\gamma$ and the exposure step of $\gamma(e/\PP)$
equals that of $\gamma$.

These definitions apply in particular to a quasicycle $\gamma$ in
$\overline{\pi/\PP}$ since the latter is a subhypergraph of
$H/\PP$. In addition, they also apply to any cycle in the graph
$(\pi/\PP)^*$, by viewing it as a quasicycle in $H/\PP$. Later in this
section, we will generalise these notions to the situation where the
plane sequence has already been completely defined.

We extend the notions of $X$-bridge and $X$-antibridge defined in
Section~\ref{sec:anti-conn} as follows: given a $\pi$-solid partition
$\RR$ of $V(H)$, $e$ is an \emph{$\RR$-bridge}
(\emph{$\RR$-antibridge}) if there is $X\in\RR$ such that $e$ is an
$X$-bridge ($X$-antibridge, respectively).

We say that a hyperedge $e$ of $H$ crossing $\RR$ is \emph{redundant}
(with respect to $\pi$ and $\RR$) if $e$ is not used by $\pi$ and $e$
is not an $\RR$-antibridge. Note that a hyperedge $e$ unused by $\pi$
is redundant if $\size e = 2$, or more generally, if each of its
vertices is in a different class of $\RR$.

Furthermore, a hyperedge $e$ of $H$ is \emph{weakly redundant} (with
respect to $\pi$ and $\RR$) if either it is redundant, or it is used
by $\pi$ and is not an $\RR$-bridge.

\begin{figure}
  \hf\sfig{27}{}\hf\sfig{28}{}\hf\\
  \hf\qquad\sfig{29}{}\hf\sfig{30}{}\hspace{1em}
  \caption{Constructing the partition $\ptn\pi 1 0$. (a) A labelling
    of the classes of $\PP := \ptn\pi0\infty$ and two hyperedges in
    the example from Figure~\ref{fig:seq}. (b) The quasigraph
    $\pi/\PP$ (bold) in $H/\PP$. (c) A quasicycle $\gamma$ (bold) in
    $\overline{\pi/\PP}$. The leading hyperedges of $\gamma$ are $e$
    and $f$. (d) The partition $\ptn\pi10$ (assuming that
    $f\leq_E e$). Note that $e$ is a $\PP$-antibridge.}
  \label{fig:limit}
\end{figure}

We are now ready to define the partition $\ptn\pi i 0$ (see
Figure~\ref{fig:limit} for an illustration). We will say that
$\ptn\pi i 0$ is obtained from $\PP$ by the \emph{limit step
  $(i-1,\infty)$}. (Recall that $\PP$ denotes the partition
$\ptn\pi{i-1}\infty$.) At the same time, we will define the
\emph{decisive hyperedge at $(i-1,\infty)$}, $d^\pi_{i-1}$, for the
current limit step. This will be a hyperedge of $H$; for technical
reasons, we also allow two extra values, $\termx$ and $\stopx$.

If the complement of $\pi/\PP$ in $H/\PP$ is acyclic, we define
$\ptn\pi i 0 = \PP$ and say that $\pi$ \emph{terminates at
  $(i-1,\infty)$}. We set $d^\pi_{i-1} = \termx$. 

Otherwise, let $L$ be the set of hyperedges $f$ of $H$ for which
there exists a quasicycle $\gamma$ in $\overline{\pi/\PP}$ such that
$f$ is a leading hyperedge of $\gamma$. We define $\ptn\pi i 0 = \PP$
if $L$ contains a weakly redundant hyperedge $f$ (with respect to
$\pi$ and $\PP$). In this case, we say that $\pi$ \emph{stops at
  $(i-1,\infty)$}; we set $d^\pi_{i-1} = \stopx$.

If no weakly redundant hyperedge exists in $L$, choose the maximum
hyperedge $e$ in $L$ according to a fixed linear ordering $\leq_E$ of
all hyperedges in $H$. This case is illustrated in
Figure~\ref{fig:limit}. (For the purposes of this and the following
section, the choice of $\leq_E$ is not important; it will be discussed
in more detail in Section~\ref{sec:bad}.) Set $d^\pi_{i-1}$ equal to
$e$. We say that $\pi$ \emph{continues} at $(i-1,\infty)$. Moreover,
in this case, any quasicycle in $\overline{\pi/\PP}$ whose leading
hyperedge is $e$ will be referred to as a \emph{decisive quasicycle at
  $(i-1,\infty)$}.

Since $e$ is not weakly redundant, we can distinguish the following
two cases for the definition of $\ptn\pi i 0$:
\begin{itemize}
\item if $e$ is a $\PP$-antibridge, then the classes of $\ptn\pi i 0$
  are all the anticomponents of $\pi$ on $X$ in $H-e$, where $X$
  ranges over all classes of $\PP$,
\item if $e$ is a $\PP$-bridge, then the classes of $\ptn\pi i 0$ are
  all the components of $\pi-e$ on $X$ in $H$, where $X$ ranges over
  all classes of $\PP$.
\end{itemize}
Note that in the first of these cases, $e$ is not used by $\pi$, while
in the second case, it is used by $\pi$.

The subsequent partitions $\ptn\pi i 1,\ptn\pi i 2,\dots$ are then
defined as described above, and the partition $\ptn\pi i\infty$ is
defined analogously to $\ptn\pi 0\infty$. Iter\-ating, we obtain the
whole plane sequence $\lseq\pi$ and the partitions
$\ptn\pi 0\infty,\ptn\pi 1 \infty,\dots$. By the finiteness of $H$,
there is some $i_0$ such that
$\ptn\pi {i_0} \infty = \ptn\pi {i_0+1} \infty$, and we define
$\ptn\pi\infty\infty$ as $\ptn\pi {i_0} \infty$.

Observe that in the cases where $\pi$ terminates or stops at
$(i-1,\infty)$ and we define $\ptn \pi i 0$ as
$\ptn \pi {i-1} \infty$, this partition will in fact equal
$\ptn\pi\infty\infty$ since none of the subsequent steps in the
construction of the plane sequence of $\pi$ will lead to any
modifications.

Now that the sequence of a quasigraph $\pi$ has been completely
defined, let us revisit the definitions of the terms `exposure step'
and `leading hyperedge'. Although these are defined relative to a
partition $\ptn\pi{i-1}\infty$ (for some $i\geq 1$), this only affects
the scope of the definitions: for instance, if vertices $u,v$ are
contained in different classes of $\ptn\pi\ell\infty$, where
$\ell\geq i$, then the exposure step of $uv$ is the same whether we
use $\ptn\pi{i-1}\infty$ or $\ptn\pi\ell\infty$ for the definition.

In particular, if we let $\QQ=\ptn\pi\infty\infty$, then it makes
sense to speak of leading hyperedges of any quasicycle in
$\overline{\pi/\QQ}$ or the exposure step of a pair of vertices
contained in different classes of $\QQ$.

We now define a partial order on the set of all quasigraphs in $H$
that is crucial for our argument. First, we define the
\emph{signature} $\sig\pi$ of a quasigraph $\pi$ as the sequence
\begin{equation*}
  \sig\pi =
  (\ptn\pi00,\ptn\pi01,\dots,\ptn\pi0\infty,d^\pi_0,\ptn\pi10,\dots,
  \ptn\pi1\infty,d^\pi_1,\dots,\ptn\pi \ell\infty,d^\pi_\ell),
\end{equation*}
where $\ell$ is minimum such that $d^\pi_\ell\in\Setx{\termx,\stopx}$.

We extend the chosen linear ordering $\leq_E$ on $E(H)$ to the union
of $E(H)$ with $\Setx{\termx,\stopx}$ by making $\termx$ the least and
$\stopx$ the greatest element, respectively. With this extension, we
are able to compare signatures in a lexicographic manner.

We derive from this an order $\worseeq$ on quasigraphs in $H$, setting
$\pi \worseeq \rho$ if $\sig\pi$ is smaller than or equal to
$\sig\rho$ in the lexicographic order on the set of signatures of
quasigraphs.

It will be convenient to define several related notions to facilitate
the comparison of quasigraphs. Let $(i,j)\in\pairs$. We define the
\emph{$(i,j)$-prefix} $\sig\pi_{(i,j)}$ of $\sig\pi$ as follows:
\begin{itemize}
\item if $j < \infty$, then $\sig\pi_{(i,j)}$ is the initial segment
  of $\sig\pi$ ending with (and including) $\ptn\pi i j$,
\item if $j = \infty$, then $\sig\pi_{(i,j)}$ is the initial segment
  of $\sig\pi$ ending with (and including) $d^\pi_i$.
\end{itemize}

We let $\pi\worseeqx i j \rho$ if $\sig\pi_{(i,j)}$ is
lexicographically smaller or equal to $\sig\rho_{(i,j)}$. Furthermore,
we define
\begin{align*}
  \pi \eq \rho &\text{\quad if\quad} \pi\worseeq\rho \text{ and
  }\rho\worseeq\pi,\\
  \pi \eqx i j \rho &\text{\quad if\quad} \pi\worseeqx i j \rho \text{ and
  }\rho\worseeqx i j \pi.
\end{align*}

Lastly, the notation $\pi\worse\rho$ means $\pi\worseeq\rho$ and
$\pi\not\eq\rho$.


\section{The main result: a variant of the Skeletal Lemma}
\label{sec:skeletal}

In this section, we are finally in a position to state and prove the
main result of this paper that is essentially a more specific version
of~\cite[Lemma 17]{KV-hamilton}. Before we state it, we need one more
definition.

Let $H$ be a $3$-hypergraph and $\pi$ an acyclic quasigraph in $H$. A
partition $\PP$ of $V(H)$ is \emph{$\pi$-skeletal} if both of the
following conditions hold:
\begin{enumerate}[(1)]
\item for each $X\in\PP$, $\pi$ is both connected on $X$ and
  anticonnected on $X$ (i.e., $\PP$ is $\pi$-solid),
  \item the complement of $\pi/\PP$ in $H/\PP$ is acyclic.
\end{enumerate}

\begin{theorem}[Skeletal Lemma, stronger version]\label{t:enhancing}
  Let $\pi$ be a quasigraph in a 3-hypergraph $H$. If
  $\ptn\pi\infty\infty$ is not $\pi$-skeletal or $\pi$ is not acyclic,
  then there is a quasigraph $\rho$ in $H$ such that either
  $\pi \sqsubset \rho$, or $\rho\eq\pi$ and $\rho$ uses fewer
  hyperedges than $\pi$.
\end{theorem}

An obvious corollary of Theorem~\ref{t:enhancing} (which will be
further strengthened in Section~\ref{sec:bad}) is the following:
\begin{corollary}\label{cor:main}
  For any $3$-hypergraph $H$, there exists an acyclic quasigraph $\pi$
  such that $\ptn\pi\infty\infty$ is $\pi$-skeletal.
\end{corollary}

Before proving Theorem~\ref{t:enhancing}, we need to establish the
following crucial lemma. The situation is illustrated in
Figure~\ref{fig:qc-addition}.

\begin{figure}
  \centering\fig{11}
  \caption{The situation in Lemma~\ref{l:qc-addition}: the quasigraph
    $\pi$ (bold) and the partition $\QQ$ (dark gray). Only some of
    the hyperedges and vertices are shown; in particular, $\QQ$ is
    assumed to be $\pi$-solid.}
  \label{fig:qc-addition}
\end{figure}

\begin{lemma}\label{l:qc-addition}
  Let $\pi$ be a quasigraph in a 3-hypergraph $H$ and let $\QQ$ be a
  $\pi$-solid partition of $V(H)$. Suppose that $X$ is a subset of
  $V(H)$ such that $\pi$ is anticonnected on $X$ and $\QQ$ refines
  $\Setx{X,V(H)-X}$. Suppose further that $\gamma$ is a quasicycle in
  $\overline{\pi/\QQ}$ all of whose vertices are subsets of $X$ (as
  classes of $\QQ$).

  If $\gamma$ has a redundant leading hyperedge $e$ (with respect to
  $\pi$ and $\QQ$), then there are vertices $u,v\in e$ such that each
  of $u$ and $v$ is contained in a different class in $\gamma(e)$, and
  the quasigraph $\pi + (uv)_e$ is anticonnected on $X$.
\end{lemma}
\begin{proof}
  Let the vertices of $\gamma^*$ be $Q_1,\dots,Q_k\in\QQ$ in order,
  such that $\gamma(e) = \Setx{Q_k,Q_1}$.
  
  \begin{claim}\label{cl:S}
    The quasigraph $\pi$ is anticonnected on $Q_1\cup\dots\cup Q_k$
    in $H-e$.
  \end{claim}
  \begin{claimproof}
    Since $e$ is redundant, it is not a $\QQ$-antibridge, so $\pi$
    is anticonnected on each $Q_i$ in $H-e$, where $i=1,\dots,k$. We
    prove, by induction on $j$, that $\pi$ is anticonnected on
    $Q_1\cup\dots\cup Q_j$ in $H-e$, where $1\leq j\leq k$. The case
    $j=1$ is clear. Supposing that $j > 1$ and the statement is
    valid for $j-1$, we prove it for $j$. Consider two consecutive
    vertices $Q_{j-1}, Q_j$ of $\gamma^*$ and the edge $f$ of
    $\gamma^*$ joining them. Since $\gamma$ is a quasigraph in
    $\overline{\pi/\QQ}$, $f$ corresponds to a hyperedge $h\neq e$
    of $H$ intersecting both $Q_{j-1}$ and $Q_j$, and such that
    $\pi(h)$ is contained in $Q_{j-1}$ or $Q_j$ (including the case
    $\pi(h)=\emptyset$). By the induction hypothesis and
    Lemma~\ref{l:union}, $\pi$ is anticonnected on
    $(Q_1\cup\dots\cup Q_{j-1})\cup Q_j$ in $H-e$.
  \end{claimproof}
  
  Let $u\in Q_k\cap e$ and $v\in Q_1\cap e$. Observe that $u$ and
  $v$ are contained in different classes of $\gamma(e)$ as stated in
  the lemma. By Claim~\ref{cl:S}, $u$ and $v$ are contained in the
  same anticomponent of $\pi$ on $X$ in $H-e$. It follows that $u$
  and $v$ are contained in the same anticomponent $A$ of
  $\rho := \pi+(uv)_e$ on $X$ in $H$. In fact, the following holds:
  
  \begin{claim}\label{cl:rho-anti}
    The quasigraph $\rho$ is anticonnected on $X$ in $H$.
  \end{claim}
  \begin{claimproof}
    Suppose the contrary and consider a partition $\RR$ of $X$ such
    that for each hyperedge $f$ of $H$ crossing $\RR$, $\rho(f)$ is
    not contained in any class of $\RR$. Then $e$ must cross $\RR$,
    since otherwise $\RR$ would demonstrate that $\pi$ is not
    anticonnected on $X$, contrary to the assumption of the
    lemma. Thus, $\rho(e)$ is not contained in any class of $\RR$,
    and since $u,v\in X$, $u$ and $v$ are contained in distinct
    classes of $\RR$. The partition $\RR[A]$ of $A$ (the
    anticomponent of $\rho$ defined above) induced by $\RR$ is
    therefore nontrivial. Since $\rho$ is anticonnected on $A$,
    there is a hyperedge $h$ of $H$ such that $h$ crosses $\RR[A]$
    and $\rho(h)$ is contained in a class of $\RR[A]$. But then $h$
    crosses $\RR$ while $\rho(h)$ is contained in a class of $\RR$,
    a contradiction with the choice of $\RR$ which proves the claim.
  \end{claimproof}
  
  We have shown that the present choice of $u$ and $v$ satisfies all
  requirements of the lemma. This concludes the proof.
\end{proof}

The following lemma is essential to relating the decisive quasicycles involved
in the construction of the plane sequences of quasigraphs $\pi$ and
$\pi-e$, where $e$ is a hyperedge. Note that it can equally well be
used to quasigraphs $\pi$ and $\pi + (uv)_e$, where $u,v\in e$. 

\begin{lemma}\label{l:qc}
  Suppose that $\RR$ is a partition of $V(H)$. If $e$ is a hyperedge
  of $H$ with $\pi(e)$ contained in a class of $\RR$, then
  \begin{equation*}
    \pi/\RR = (\pi-e)/\RR.
  \end{equation*}
  In particular, the complements of $\pi/\RR$ and of $(\pi-e)/\RR$
  coincide.
\end{lemma}
\begin{proof}
  Both $\pi/\RR$ and $(\pi-e)/\RR$ are quasigraphs in $H/\RR$. Their
  values on $f/\RR$, where $f\neq e$ is any $\RR$-crossing hyperedge,
  are clearly the same. It thus suffices to compare the values of
  $\pi/\RR$ and $(\pi-e)/\RR$ on $e/\RR$ under the assumption that $e$
  is $\RR$-crossing. The latter value is $\emptyset$ since $e$ is not
  used by $\pi-e$. The former one is also $\emptyset$ since $\pi(e)$
  is not $\RR$-crossing.
\end{proof}

Given $(i,j)\in\pairs$ with $i,j < \infty$ and $(i,j)\neq (0,0)$, the
\emph{predecessor} of the partition $\ptn\pi i j$ is the partition
$\ptn\pi i {j-1}$ if $j > 0$, and $\ptn\pi {i-1}\infty$ if $j = 0$ and
$i > 0$. The predecessors of the other partitions in the sequence for
$\pi$ are undefined.

\begin{observation}\label{obs:exposed}
  Let $\QQ$ be a partition of $V(H)$. If a quasicycle $\gamma$ in
  $\overline{\pi/\QQ}$ is exposed at $(i,j)$ with respect to $\pi$,
  then $i,j < \infty$ and $(i,j)\neq (0,0)$; in particular, the
  predecessor of $\ptn\pi i j$ exists. In addition, if $j = 1$ and
  $i\geq 1$, then $d^\pi_{i-1}$ is a hyperedge not used by $\pi$.
\end{observation}
\begin{proof}
  Clearly, $(i,j) \neq (0,0)$ since $\ptn\pi00 = \Setx{V(H)}$. By the
  definition of the sequence of $\pi$, for any $r\geq 0$,
  $\ptn\pi r \infty$ is equal to one (actually, infinitely many) of
  the partitions $\ptn\pi r s$, where $0\leq s < \infty$, and hence it
  cannot be the exposing partition for $\gamma$. A similar argument
  applies to $\ptn\pi\infty\infty$. As for the last statement, suppose
  that the exposing partition is $\ptn\pi i 1$. Clearly,
  $d^\pi_{i-1}\notin\Setx{\termx,\stopx}$, so it is a hyperedge. If it
  were used by $\pi$, it would be a $\ptn\pi{i-1}\infty$-bridge, and
  the classes of $\ptn\pi i 0$ would be the components of
  $\pi-d^\pi_{i-1}$ on the classes of $\ptn\pi{i-1}\infty$. By the
  definition of the sequence of $\pi$, $\ptn\pi i 1 = \ptn\pi i 0$ and
  so $\ptn\pi i 1$ cannot be the exposing partition for $\gamma$.
\end{proof}

\begin{observation}
  \label{obs:anti}
  Suppose that $\pi$ is a quasigraph in $H$, $X\subseteq V(H)$ and $e$
  is a hyperedge such that $e\cap X = \Setx{u,v}$. If $\pi+(uv)_e$ is
  anticonnected on $X$, then $e$ is not an $X$-antibridge with respect
  to $\pi$. 
\end{observation}
\begin{proof}
  Assume the contrary. Then there is a nontrivial partition $\RR$ of
  $X$ such that for each $\RR$-crossing hyperedge $f$ of $H-e$,
  $\pi(f)$ is not contained in any class of $\RR$. Since there is no
  such partition for $\pi+(uv)_e$ in $H$, it must be that $e$ crosses
  $\RR$ and $\Setx{u,v}$ is a subset of a class of $\RR$. That is
  impossible since $e\cap X = \Setx{u,v}$.
\end{proof}

\begin{lemma}\label{l:stable}
  Let $\pi$ be a quasigraph in $H$. Let $e$ be a hyperedge not used by
  $\pi$ such that vertices $u,v\in e$ are contained in different
  classes of a partition $\ptn\pi i j$ (where $0 \leq i,j < \infty$),
  but both of them are contained in the same class $X$ of its
  predecessor. If the quasigraph $\pi + (uv)_e$ is anticonnected on
  $X$, then $\pi \sqsubset \pi + (uv)_e$.
\end{lemma}
\begin{proof}
  Let $\rho = \pi + (uv)_e$. Suppose that $\pi\not\sqsubset\rho$. We
  begin by proving the following claim:
  \begin{equation}\label{eq:stable} \text{$\pi \worseeqx s t
      \rho$ for all $(s,t)\in\pairs$ with $(s,t) \leq (i,j)$.}
  \end{equation}

  We proceed by (transfinite) induction on $(s,t)$, assuming the claim
  for all smaller pairs in $\pairs$. The claim~\eqref{eq:stable} holds
  for $(s,t) = (0,0)$; assume therefore that $(s,t) > (0,0)$. Suppose
  first that $0 < t < \infty$ and the statement holds for
  $(s,t-1)$. If $t$ is odd, then any class $A$ of $\ptn\pi{s}{t}$ is a
  component of $\pi$ on a class of $\ptn\pi{s}{t-1}$. We have either
  $A\supseteq X$, or $A\cap X = \emptyset$. In both cases, $\rho$ is
  clearly connected on $A$ (by Lemma~\ref{l:sub-add}(ii) in the former
  case). Thus, $\ptn\pi {s}{t} \leq \ptn\rho{s}{t}$ and
  $\pi\worseeqx s t \rho$.

  Next, if $t$ is even (and nonzero), we proceed similarly: if $A$ is
  an anticomponent of $\pi$ on a class of $\ptn\pi{s}{t-1}$ and
  $A\supseteq X$, then $\rho$ is anticonnected on $A$ by
  Lemma~\ref{l:sub-add}(i), while if $A\cap X = \emptyset$, the same
  is true for trivial reasons. Thus, again, $\pi\worseeqx s t \rho$.

  The next case is $t = \infty$. By the induction hypothesis,
  $\ptn\pi s \infty \leq \ptn\rho s \infty$; without loss of
  generality, we may assume that
  $\ptn\pi s \infty = \ptn\rho s \infty$. Let
  $\SSS = \ptn\pi s \infty$. We need to show that
  $d^\pi_s \leq_E d^\rho_s$. By Lemma~\ref{l:qc},
  $\overline{\pi/\SSS} = \overline{\rho/\SSS}$. In particular, the two
  hypergraphs have the same quasicycles, and these quasicycles have
  the same leading hyperedges. It follows that
  $d^\pi_s \leq_E d^\rho_s$ or $d^\pi_s = \stopx$ --- but the latter
  does not hold since $\pi$ cannot stop (or terminate) at $(s,\infty)$
  as $\ptn\pi i j$ differs from its predecessor and
  $(s,\infty) < (i,j)$.

  It remains to consider the case $t=0$. Here, we have $s > 0$ and the
  induction hypothesis implies that $\pi\worseeqx
  {s-1}\infty\rho$. Let us assume that $\pi\eqx {s-1}\infty\rho$ and
  denote $\ptn\pi{s-1}\infty$ by $\SSS$. We want to show that $\ptn\pi
  s 0 \leq \ptn\rho s 0$.

  Let $f := d^\pi_{s-1}$. For the same reason as above, $\pi$
  continues at $(s-1,\infty)$, so $f\notin\Setx{\termx,\stopx}$. This
  means that $f$ is an $\SSS$-bridge or an $\SSS$-antibridge with
  respect to $\pi$. In addition, $d^\rho_{s-1} = f$ by the assumption
  that $\pi\eqx {s-1}\infty\rho$. Hence, $f$ is an $\SSS$-bridge or an
  $\SSS$-antibridge with respect to $\rho$ as well.
  
  Consider the set $X$ from the statement of the lemma. Let $S$ be the
  class of $\SSS$ containing $X$. Since the anticomponents of $\pi-f$
  and $\rho-f$ on $S$ are clearly the same if $\size{f\cap S}\leq 1$,
  we may assume that $\size{f\cap S}\geq 2$ --- indeed, since $f$
  crosses $\SSS$, we must have equality.

  We distinguish three cases:
  \begin{enumerate}[(a)]
  \item $f = e$,
  \item $f\neq e$ and $f$ is not used by $\pi$,
  \item $f\neq e$ and $f$ is used by $\pi$.
  \end{enumerate}

  In case (a), $f$ ($=e$) is not used by $\pi$, so it is an
  $S$-antibridge with respect to $\pi$. Since it intersects $S$ in two
  vertices, these two vertices are $u$ and $v$, and each of them is
  contained in a different anticomponent of $\pi$ on $S$ in
  $H-f$. Moreover, it must be that $i = s$, $j = 0$ and $X = S$ (since
  $u,v$ are contained in distinct classes of $\ptn\pi s 0$ but in one
  class $S$ of its predecessor). However, an assumption of the lemma
  is that $\rho$ is anticonnected on $X$, a contradiction with
  Observation~\ref{obs:anti}. In other words, case (a) cannot occur.
  
  In case (b), which is illustrated in Figure~\ref{fig:case-b}, $f$ is
  also an $S$-antibridge with respect to $\pi$ (and $\rho$). To prove
  that $\ptn\pi s 0 \leq \ptn\rho s 0$, it is enough to show that
  $\rho$ is anticonnected on each anticomponent of $\pi$ on $S$ in
  $H-f$. Let $A$ be such an anticomponent. We may assume that
  $\size{e\cap A} \geq 2$, otherwise $\rho$ is clearly anticonnected on
  $A$. Thus, $X\subseteq A$. By Lemma~\ref{l:sub-add}(i), $\rho$ is
  anticonnected on $A$ as claimed. The discussion of case (b) is
  complete.

  \begin{figure}
    \centering\fig{16}
    \caption{Case (b) in the proof of Lemma~\ref{l:stable}. A class
      $S$ of the partition $\SSS$ is shown in lightest
      gray. Hyperedges $e$ and $f$ are shown dotted. Medium gray
      regions represent the anticomponents of $\pi$ on $S$ in
      $H-f$. The set $X$ (darkest gray) is contained in an
      anticomponent $A$.}
    \label{fig:case-b}
  \end{figure}

  Lastly, in case (c), $f$ is an $S$-bridge with respect to $\pi$ and
  $\rho$. Trivially, the quasigraph $\rho-f$ is connected on each
  component of $\pi-f$ on $S$ in $H$, so
  $\ptn\pi s 0 \leq \ptn\rho s 0$.

  To summarise, each of the cases (a)--(c) leads either to a
  contradiction, or to the sought conclusion
  $\ptn\pi s 0 \leq \ptn\rho s 0$. This concludes the proof
  of~\eqref{eq:stable}.
  
  It remains to show that $\pi\sqsubset\rho$. Since $j < \infty$,
  there are three cases to distinguish based on the value of $j$.  If
  $j$ is odd, then the classes of $\ptn\pi i j\setminus\ptn\pi i {j-1}$ are
  the components of $\pi$ on $X$. Since $u$ and $v$ are in different
  classes of $\ptn\pi i j$, the replacement of $\pi$ with $\rho$ has
  the effect of adding the edge $uv$ to $\pi^*$, joining the two
  components into one. Therefore, $\ptn\pi i j < \ptn\rho i j$, and
  by~\eqref{eq:stable}, $\pi\sqsubset\rho$.

  If $j$ is even and $j > 0$, the classes of $\ptn\pi i j$ are the
  anticomponents of $\pi$ on $X$. Let $A_1$ and $A_2$ be such
  anticomponents containing $u$ and $v$, respectively. By
  Lemma~\ref{l:union}, $\pi$ is anticonnected on $A_1\cup A_2$ since
  $e$ intersects both of these sets and is not used by $\pi$. This
  contradiction means that the present case is not possible.

  Finally, if $j=0$, $uv$ is exposed at $(i,0)$ with respect to
  $\pi$. Let $f = d^\pi_{i-1}$ be the corresponding decisive
  hyperedge. Since $\pi$ continues at $(i-1,\infty)$, $f$ is an
  $X$-bridge or an $X$-antibridge with respect to $\pi$. This shows
  that $f\neq e$ because $e$ is not an $X$-antibridge by
  Observation~\ref{obs:anti}, and is not an $X$-bridge because it is
  not used by $\pi$. Furthermore, similarly to the preceding case, it
  cannot be that $f$ is an $X$-antibridge, for then $e$ would
  intersect two anticomponents of $\pi$ on $X$ in $H-f$, which is
  impossible by Lemma~\ref{l:union}.

  Thus, $f$ is an $X$-bridge with respect to $\pi$, and the classes of
  $\ptn\pi i 0$ are the components of $\pi-f$ on $X$ in $H$. Since
  $u,v$ are in different components, the quasigraph $\rho-f$ is
  connected on $X$ in $H$, and consequently $\rho$ stops at
  $\ptn\rho{i-1}\infty$ and $\pi\sqsubset\rho$. This concludes the
  proof.
\end{proof}

\begin{lemma}\label{l:stable-cycle}
  Let $\pi$ be a quasigraph in $H$. Let $e$ be a hyperedge of $H$ used
  by $\pi$ and $X\subseteq V(H)$ such that $\pi(e)\subseteq X$, $\pi-e$ is
  connected on $X$, and one of the following conditions is satisfied:
  \begin{enumerate}[(a)]
  \item the vertices of $\pi(e)$ are contained in different classes of
    $\ptn\pi i j$ $(0\leq i,j < \infty)$ and $X$ is a class of its
    predecessor,
  \item $X\in\ptn\pi\infty\infty$ and $\ptn\pi\infty\infty$ is
    $\pi$-skeletal.
  \end{enumerate}
  Then $\pi \sqsubseteq \pi-e$. In addition, if (a) is satisfied, then
  $\pi \sqsubset \pi-e$.
\end{lemma}
\begin{proof}
  Suppose that condition (a) is satisfied. We prove, by (transfinite)
  induction on $(s,t)$, that
  \begin{equation}
    \label{eq:1}
    \pi \worseeqx s t \pi-e \text{ for all $(s,t)\in\pairs$ with $(s,t) \leq (i,j)$}.
  \end{equation}
  The statement holds if $(s,t)=(0,0)$. Consider $(s,t) > (0,0)$. If
  $0 < t < \infty$, then we may suppose that $\pi \eqx s {t-1} \pi-e$; by
  Lemma~\ref{l:sub-remove}, $\ptn\pi s t \leq \ptn{\pi-e} s t$ and
  therefore $\pi \worseeqx s t \pi-e$ as desired.

  Suppose that $t = \infty$. Let $\PP = \ptn\pi s\infty$. Without loss
  of generality, $\PP = \ptn{\pi-e} s \infty$. Since $X$ is a subset
  of a class of $\PP$, Lemma~\ref{l:qc} implies that the complement of
  $\pi/\PP$ is the same as the complement of $(\pi-e)/\PP$. In
  particular, the quasicycles in these hypergraphs, as well as the
  sets of their leading hyperedges, are the same.

  We state the following simple observation as a claim for easier
  reference later in the proof:
  \setcounter{claim}0
  \begin{claim}\label{cl:no-wr}
    No leading hyperedge of any quasicycle in the complement of
    $\pi/\PP$ is weakly redundant.
  \end{claim}
  Indeed, if the claim did not hold, then $\pi$ would have stopped at
  $(s,\infty)$, but condition (a) implies that $\ptn\pi i j$ differs
  from its predecessor. Since $(s,\infty) < (i,j)$, this would be a
  contradiction.

  \medskip
  By Claim~\ref{cl:no-wr}, if $\pi-e$ stops at $(s,\infty)$, then
  $\pi\worse\pi-e$. We may therefore assume that $\pi-e$ continues at
  $(s,\infty)$, in which case the decisive hyperedges at $(s,\infty)$
  for $\pi$ and $\pi-e$ coincide. We conclude that
  $\pi\worseeqx s \infty\pi-e$ in this case.

  It remains to consider the case $t=0$. It suffices to show that
  $\ptn\pi s 0 \leq \ptn{\pi-e} s 0$ assuming that
  $\pi\eqx{s-1}\infty\pi-e$. Let $f := d^\pi_{s-1}$ be the decisive
  hyperedge at $(s-1,\infty)$ with respect to $\pi$. The assumption
  implies that $f = d^{\pi-e}_{s-1}$. Moreover,
  $f\notin\Setx{\termx,\stopx}$ because $\ptn\pi i j$ differs from its
  predecessor, so $\pi$ has to continue at $(s-1,\infty)$ as
  $(i,j) \geq (s,0)$.

  If there is a class $A\in\ptn\pi s 0$ with $\pi(e)\subseteq A$, then
  we have $X\subseteq A$ (where $X$ is the set from the lemma). It
  follows from Lemma~\ref{l:sub-remove} that $\pi-e$ is connected on
  $A$ whenever $\pi$ is, and similarly for
  anticonnectivity. Consequently, $\ptn\pi s 0 \leq \ptn{\pi-e} s 0$.

  We may therefore assume that
  \begin{equation}
    \label{eq:2}
    \pi(e) \text{ intersects two classes $Y,Y'$ of $\ptn\pi s 0$}.
  \end{equation}
  The situation is illustrated in
  Figure~\ref{fig:0}. By condition (a), $(i,j)=(s,0)$ and
  $X\in\ptn\pi{s-1}\infty$. 

  By the construction of $\ptn\pi s 0$, there are two possibilities:
  either (i) $f$ is an $X$-antibridge with respect to $\pi$ and $Y,Y'$
  are the anticomponents of $\pi$ on $X$ in $H-f$, or (ii) $f$ is an
  $X$-bridge with respect to $\pi$ and $Y,Y'$ are the components of
  $\pi-f$ on $X$ in $H$.

  Note first that $f\neq e$: otherwise, $f$ would necessarily be an
  $X$-bridge (since $e$ is used by $\pi$), but $\pi-e$ is assumed to
  be connected on $X$.

  We next rule out possibility (ii). Assume that $Y,Y'$ are the
  components of $\pi-f$ on $X$ in $H$. Now $\pi(e)$ intersects $Y$ and
  $Y'$, so $\pi-f$ is connected on $Y\cup Y'$ thanks to $e$, a
  contradiction.

  Let us therefore restrict ourselves to possibility (i). Let $\gamma$
  denote a quasicycle in the complement of $\pi/\ptn\pi{s-1}\infty$
  such that $f$ is a leading hyperedge of $\gamma$. By
  Lemma~\ref{l:qc}, $\gamma$ is a quasicycle in the complement of
  $(\pi-e)/\ptn\pi{s-1}\infty$ as well.

  Since $e$ intersects both $Y$ and $Y'$, (i) implies that $\pi-e$ is
  anticonnected on $Y\cup Y'$ in $H-f$. Consequently, the leading
  hyperedge $f$ of $\gamma$ is \emph{redundant} with respect to
  $\pi-e$ and $\ptn\pi{s-1}\infty$. This contradicts our assumption
  that $\pi\eqx{s-1}\infty\pi-e$ since $\pi$ continues at
  $(s-1,\infty)$ while $\pi-e$ stops there. The contradiction means
  that if $\pi(e)$ intersects two classes of $\ptn\pi s 0$, then
  actually $\pi\worsex{s-1}\infty\pi-e$. The case $t=0$ is settled.

  \begin{figure}
    \centering\fig{17}
    \caption{A situation in the proof of Lemma~\ref{l:stable-cycle}
      when $\pi(e)$ intersects two classes $Y,Y'$ of $\ptn\pi s0$
      (shown in darker gray). Light gray represents a class $X$ of
      $\ptn\pi{s-1}\infty$. Observe that $\pi-e$ is anticonnected on
      $Y\cup Y'$ in $H-f$.}
    \label{fig:0}
  \end{figure}
  
  Having proved~\eqref{eq:1}, let us now show that
  \begin{equation*}
    \pi\worsex i j\pi-e,
  \end{equation*}
  where $i,j$ are as in the statement of the lemma. By assumption,
  $\pi(e)$ intersects two classes of $\ptn\pi i j$. We note that
  $0\leq j < \infty$ and consider two possibilities: $j = 0$ and
  $j > 0$ ($j$ finite). If $j = 0$, then we have seen in the above
  paragraph that if $\pi(e)$ intersects two classes of $\ptn\pi s 0$
  (for any $s\leq i$), then $s=i$ and $\pi\worsex{s-1}\infty\pi-e$. In
  particular, $\pi\worsex i j\pi-e$.

  Suppose now that $j > 0$ is finite. Without loss of generality,
  $\pi\eqx i {j-1}\pi-e$. Since $\pi(e)$ intersects two classes of
  $\ptn\pi i j$, we see from the definition of the sequence of $\pi$
  that $j$ is even. Hence, the vertices of $\pi(e)$ lie in different
  anticomponents of $\pi$ on $X$, and $\pi-e$ is anticonnected on the
  union of these anticomponents by Lemma~\ref{l:union}. Consequently,
  $\ptn\pi i j < \ptn{\pi-e} i j$ and $\pi\worsex i j \pi-e$. This
  concludes the proof for condition (a).
  
  A very similar inductive proof to that used to prove~\eqref{eq:1}
  works when condition (b) is satisfied. The main difference is that
  in the $t = \infty$ case, Claim~\ref{cl:no-wr} now holds for a
  different reason, namely that $\ptn\pi\infty\infty$ is
  $\pi$-skeletal, which means that $\pi$ never stops --- consequently,
  there is no weakly redundant leading hyperedge of a quasicycle in
  $\overline{\pi/\ptn\pi{s-1}\infty}$.

  Another difference is that condition~\eqref{eq:2} cannot occur if
  (b) is satisfied, which makes the proof for condition (b) somewhat
  shorter.
\end{proof}
  
We can now proceed to the proof of Theorem~\ref{t:enhancing}.

\begin{proof}[\textbf{Proof of Theorem~\ref{t:enhancing}}]
  Let $\QQ=\ptn\pi\infty\infty$. We distinguish the following cases.
  
  \begin{xcase}{$\QQ$ is not $\pi$-skeletal.}%
    By the construction, it is clear that $\QQ$ is $\pi$-solid. Thus,
    $\overline{\pi/\QQ}$ contains a quasicycle. Consider the least $s$
    such that $\ptn\pi s \infty = \QQ$. Since $\pi$ stops at
    $\ptn\pi s \infty$, there is a quasicycle $\gamma$ in
    $\overline{\pi/\QQ}$ and a leading hyperedge $e$ of $\gamma$ such
    that $e$ is weakly redundant. Let the exposure step for $\gamma$
    be $(i,j)$, where $0\leq i,j < \infty$. We put
    $\PP = \ptn\pi i j$, and let $X$ be the class of the predecessor
    of $\ptn\pi i j$ containing both vertices of
    $\gamma(e)$. Furthermore, let $Q_1,Q_2\in\QQ$ be the two vertices
    of $\gamma(e)$, and let $P_i\in\PP$ be such that
    $P_i\supseteq Q_i$ ($i = 1,2$).

    \begin{xsubcase}{$e$ is not used by $\pi$.}%
      First, $j$ is odd or $j = 0$; otherwise, $P_1,P_2$ would be
      anticomponents of $\pi$ on $X$, but Lemma~\ref{l:union} shows
      that $\pi$ is anticonnected on $P_1\cup P_2$, which would be a
      contradiction.

      Consider first the case that $j$ is odd, so $P_1$ and $P_2$ are
      components of $\pi$ on $X$. Moreover, suppose for now that
      $j > 1$. By the construction of the sequence for $\pi$, $\pi$ is
      anticonnected on $X$. Applying Lemma~\ref{l:qc-addition} (with
      the current values of $X$, $Q$, $\gamma$ and $e$), we obtain
      vertices $u,v$ such that $u\in Q_1$, $v\in Q_2$ and
      $\pi + (uv)_e$ is anticonnected on $X$. Lemma~\ref{l:stable}
      then implies that $\pi \sqsubset \pi + (u_1u_2)_e$ and we are
      done.

      Suppose that $j = 1$. If $i=0$ or the decisive hyperedge
      $d^\pi_{i-1}$ at $(i-1,\infty)$ is not used by $\pi$, then the
      above argument works, since in this case $\pi$ is anticonnected
      on $X$. The other case ($i > 0$ and $d^\pi_{i-1}$ is used by $\pi$) is
      excluded by Observation~\ref{obs:exposed}.
      This settles the case $j=1$ and more broadly the case that $j$
      is odd.

      It remains to consider the possibility that $j = 0$. Clearly,
      $i > 0$ as $\ptn\pi00 = \Setx{V(H)}$. The predecessor of
      $\ptn\pi i 0$ is $\ptn\pi{i-1}\infty$, which is
      $\pi$-solid. Thus, $\pi$ is anticonnected on $X$ and the
      argument used for odd $j>1$ applies.%
    \end{xsubcase}
    \begin{xsubcase}{$e$ is used by $\pi$.}%
      Since $e$ is a leading hyperedge of a quasicycle in
      $\overline{\pi/\QQ}$, $\pi(e)$ is a subset of some
      $Q\in\QQ$. Since $e$ is weakly redundant with respect to $\QQ$,
      it is not a $Q$-bridge --- in other words, $\pi-e$ is connected
      on $Q$. Lemma~\ref{l:stable-cycle} implies that
      $\pi \sqsubseteq \pi-e$. Since $\pi-e$ uses fewer hyperedges
      than $\pi$, it has the desired properties.
    \end{xsubcase}%
  \end{xcase}

  \begin{xcase}{$\pi$ is not acyclic.}%
    Suppose that there is a cycle $C$ in the graph
    $\pi^*$. That means that either $\pi/\QQ$ is not acyclic, or there
    is a cycle in the induced subgraph of $\pi^*$ on some $Q\in\QQ$.

    \begin{xsubcase}{The quasigraph $\pi/\QQ$ is not acyclic.}%
      Let $C$ be a cycle in $(\pi/\QQ)^*$ and suppose its exposure
      step is $(i,j)$, where $0 \leq i,j < \infty$. Let $e$ be a
      leading hyperedge of $C$. Note that each of the vertices of
      $\pi(e)$ is in a different class of the partition $\ptn\pi i j$,
      but both are contained in the same class $X$ of its predecessor.

      Since $e$ is contained in the cycle $C$ all of whose vertices
      are subsets of $X$, $\pi-e$ is connected on $X$. Thus, the
      assumptions of Lemma~\ref{l:stable-cycle} are satisfied (we use
      the current values of $e$ and $X$). It follows that
      $\pi \sqsubseteq \pi-e$; since $\pi-e$ uses fewer hyperedges than $\pi$,
      $\pi-e$ has the desired properties.
    \end{xsubcase}

    \begin{xsubcase}{There is $Q\in\QQ$ such that the induced subgraph
      of $\pi^*$ on $Q$ contains a cycle.}%
    Let $C$ be a cycle in the induced subgraph of $\pi^*$ on $Q$, and
    let $e$ be a hyperedge such that $\pi(e)$ is an edge of
    $C$. Clearly, $\pi-e$ is connected on $Q$. By
    Lemma~\ref{l:stable-cycle}, $\pi \sqsubseteq \pi-e$. Since $\pi-e$
    uses fewer hyperedges than $\pi$ does, we are done.%
    \end{xsubcase}%
  \end{xcase}%
  \vspace{-1.2em}
\end{proof}


\section{Removing bad leaves}
\label{sec:bad}

For the purposes of the application of the Skeletal Lemma
in~\cite{KV-ess-9}, we have to prove the lemma in a stronger form
(Theorem~\ref{t:enhancing}) allowing us to deal with a certain
configuration that is problematic for the analysis in~\cite{KV-ess-9},
namely a `bad leaf' in a quasigraph. As a result, we will be able to
exclude this configuration in Theorem~\ref{t:no-bad} (at the cost of
some local modifications to the hypergraph).

Let $H$ be a 3-hypergraph and let $\pi$ be an acyclic quasigraph in
$H$. In each component of the graph $\pi^*$, we choose an arbitrary
root and orient all the edges of $\pi^*$ toward the root. A hyperedge
$e$ of $H$ is \emph{associated with} a vertex $u$ if it is used by
$\pi$ and $u$ is the tail of $\pi(e)$ in the resulting oriented
graph. Thus, every vertex has at most one associated hyperedge, and
conversely, each hyperedge is associated with at most one vertex.

A vertex $u$ of $H$ is a \emph{bad leaf} for $\pi$ if all of the
following hold:
\begin{enumerate}[\quad(i)]
\item $u$ is a leaf of $\pi^*$,
\item $u$ is incident with exactly three hyperedges, exactly one of
  which has size 3 (say, $e$), and
\item $e$ is associated with $u$.
\end{enumerate}

\begin{figure}
  \centering
  \hf
  \sfig1{}\hf\sfig2{}
  \hf
  \caption{(a) A bad leaf $u$. (b) The result of the switch at $u$. }
  \label{fig:bad}
\end{figure}

To eliminate a bad leaf $u$, we use a \emph{switch} operation
illustrated in Figure~\ref{fig:bad}. Suppose that $u$ is incident
with hyperedges $ua$, $ub$ and $ucd$, where $ucd$ is associated with
$u$, and $\pi(ucd) = uc$. We remove from $H$ the hyperedges $ua$, $ub$
and $ucd$ and add the hyperedges $uab$, $uc$ and $ud$; the resulting
hypergraph is denoted by $H^{(u)}$. We say that a hypergraph $\tilde
H$ is \emph{related} to $H$ if it can be obtained from $H$ by a finite
series of switch operations.

With $\pi$ as above, a quasigraph $\pi^{(u)}$ in $H^{(u)}$ is obtained
by setting $\pi^{(u)}(uc) = uc$, and leaving both $ud$ and $uab$
unused. Observe that $(\pi^{(u)})^* = \pi^*$, and $\pi^{(u)}$ has
fewer bad leaves than $\pi$.

A problem we have to address is that a partition $\PP$ which is
$\pi$-skeletal in $H$ need no longer be $\pi^{(u)}$-skeletal in
$H^{(u)}$, since the switch may create an unwanted cycle in
$\overline{\pi^{(u)}/\PP}$. This is illustrated in
Figure~\ref{fig:switch-cycle}. The following paragraphs describe the
steps taken to resolve this problem. First, we extend the order
$\sqsubseteq$ defined on quasigraphs in $H$ to the set of all
quasigraphs in hypergraphs related to $H$. Since all such hypergraphs
have the same vertex set, we can readily compare the partitions of
their vertex sets.

\begin{figure}
  \centering\hf\subfloat[]{\fig{12}}\hf\subfloat[]{\fig{13}}\hf
  \caption{(a) A quasigraph $\pi$ in $H$ such that the complement of
    $\pi$ is acyclic. (b) The quasigraph obtained by a switch at the
    vertex $u$ no longer has acyclic complement.}
  \label{fig:switch-cycle}
\end{figure}

We have to be more careful, however, in the definition of the sequence
of $\pi$, where a linear ordering $\leq_E$ of hyperedges of $H$ is
used: this ordering should involve all hyperedges of hypergraphs
related to $H$. We define $\leq_E$ as follows. We fix a linear
ordering $\leq$ of $V(H)$. On the set of 3-hyperedges of hypergraphs
related to $H$, $\leq_E$ is the associated lexicographic ordering, and
the same holds for the set of 2-hyperedges of hypergraphs related to
$H$. We make each 2-hyperedge greater than any 3-hyperedge with
respect to $\leq_E$. Finally, we add the elements $\termx$ and
$\stopx$ to the ordered set as the least and greatest element of
$\leq_E$, respectively.

This allows for a definition of the sequence of $\pi$ consistent with
the switch operation. Furthermore, the definition of the
ordering $\sqsubseteq$ as given in Section~\ref{sec:sequence} is well
suited for our purpose, and remains without change.

Let us mention that although the incorporation of the decisive hyperedges in
the signature of a quasigraph may have seemed unnecessary (mainly
thanks to Lemma~\ref{l:qc}), the present section is the reason why we
chose this definition. In fact, the only situation in our arguments
when the comparison of the decisive hyperedges is relevant is
immediately after a switch, as in Lemma~\ref{l:stable-bad} below.

We first prove that switching a bad leaf of $\pi$ does not affect the
(anti)connecti\-vi\-ty of $\pi$ on a set of vertices.

\begin{lemma}
  \label{l:sub-switch}
  Let $\pi$ be an acyclic quasigraph in $H$ and $X\subseteq
  V(H)$. Suppose that $\pi$ has a bad leaf $u$ and $\sigma$ is
  obtained from $\pi$ by switching at $u$. The following holds:
  \begin{enumerate}[\quad(i)]
  \item if $\pi$ is anticonnected on $X$, then so is $\sigma$,
  \item if $\pi$ is connected on $X$, then so is $\sigma$.
  \end{enumerate}
\end{lemma}
\begin{proof}
  We prove (i). Suppose that $\pi$ is anticonnected on $X$, but the
  quasigraph $\sigma$ (in a hypergraph $\tilde H$ related to $H$) is
  not. The definition implies that there is a partition $\PP$ of $X$
  such that for every $\PP$-crossing hyperedge $f$ of $\tilde H$,
  $\sigma(f)$ is not contained in any class of $\PP$. At the same
  time, there is a $\PP$-crossing hyperedge $e$ of $H$ such that
  $\pi(e)$ is contained in some class of $\PP$. Clearly, $e$ must be
  incident with $u$ (since the other hyperedges exist both in $H$ and
  $\tilde H$, and the values of $\pi$ and $\sigma$ coincide).

  Let the neighbours of $u$ in $\tilde H$ be labelled as in
  Figure~\ref{fig:bad}(b). Let $A$ be the class of $\PP$ containing
  $u$; by the above property of $\sigma$, we can easily see that
  $a\in A$ if $a\in X$, and similarly for $b$ and $d$. This implies
  that $c\in X\setminus A$ and $e=ucd$, but then $\pi(ucd)$ crosses
  $\PP$, a contradiction.

  Part (ii) is immediate from the fact that $\pi^* = \sigma^*$. 
\end{proof}

\begin{lemma}\label{l:stable-bad}
  Let $\pi$ be an acyclic quasigraph in $H$ such that
  $\ptn\pi\infty\infty$ is $\pi$-skeletal. If $\pi$ has a bad leaf $u$
  and the quasigraph $\sigma$ (in a hypergraph related to $H$) is
  obtained from $\pi$ by a switch at $u$, then
  $\pi \sqsubseteq \sigma$.
\end{lemma}
\begin{proof}
  We show that
  \begin{equation}
    \label{eq:stable-bad}
    \text{$\pi \worseeqx i j
      \sigma$ for all $(i,j)\in\pairs$.}
  \end{equation}
  We proceed by transfinite induction on $(i,j)$. The claim is trivial
  for $(i,j)=(0,0)$. Suppose now that $j > 0$ is finite. We may assume
  that $\pi\eqx i {j-1}\sigma$ for otherwise we are done. If $j$ is
  odd, then the classes of $\ptn\pi i j$ are the components of $\pi$
  on classes of $\ptn\pi i {j-1}$. Let $X\in\ptn\pi i {j-1}$ and let
  $A$ be a component of $\sigma$ on $X$. By
  Lemma~\ref{l:sub-switch}(ii), $\sigma$ is connected on $A$. Hence,
  $\ptn\pi i j\leq \ptn\sigma i j$, and~\eqref{eq:stable-bad}
  follows. An analogous argument, using Lemma~\ref{l:sub-switch}(i),
  can be used for even $j > 0$.

  Let us consider the case $j = \infty$. We may assume that $i$ is
  finite (otherwise, the claim follows directly from the induction
  hypothesis) and that $\pi\eqx i \infty\sigma$. Since
  $\ptn\pi\infty\infty$ is assumed to be $\pi$-skeletal, $\pi$ does
  not stop at $(i,\infty)$. Furthermore, we may assume that $\pi$ does
  not terminate at $(i,\infty)$, for otherwise we immediately conclude
  $\pi\worseeqx i \infty\sigma$ (recalling that $\termx$ is the least
  element of the ordering $\leq_E$).

  Let $\PP=\ptn\pi i\infty$ and let $\gamma$ be a quasicycle in the
  complement of $\pi/\PP$ in $H/\PP$. Define a quasigraph $\gamma'$ in
  the complement of $\sigma/\PP$ in $\tilde H/\PP$ as follows (see
  Figure~\ref{fig:switch-cases} for an illustration of several of the
  cases):
  \begin{itemize}
  \item if $\gamma$ uses a hyperedge $f/\PP$, where $f$ is a hyperedge
    of $H$ not incident with $u$, then set $\gamma'(f/\PP) =
    \gamma(f/\PP)$,
  \item if $\gamma$ uses $au/\PP$ and $bu/\PP$, then set
    $\gamma'(abu/\PP) = ab/\PP$,
  \item if $\gamma$ uses $au/\PP$ but not $bu/\PP$, then set
    $\gamma'(abu/\PP) = au/\PP$ (and symmetrically with $au$ and $bu$
    reversed),
  \item if $\gamma$ uses $ucd/\PP$ (so $u$ and $d$ are in different
    classes of $\PP$), then set $\gamma'(ud/\PP) = ud/\PP$.
  \end{itemize}

  \begin{figure}
    \hf\subfloat[]{\fig{14}}\hf\subfloat[]{\fig{15}}\hf
    \caption{Corresponding quasicycles in the complement of $\pi/\PP$
      (in $H/\PP$) and in the complement of $\pi^{(u)}/\PP$ (in
      $\tilde H/\PP$). The quasigraph $\pi$ is represented by bold
      lines, the partition $\PP$ is shown in gray. (a) A quasicycle
      $\gamma$ in the complement of $\pi/\PP$ (shown with darker gray
      bars). (b) The corresponding quasicycle $\gamma'$ in
      $\pi^{(u)}/\PP$ (shown in the same way).}
    \label{fig:switch-cases}
  \end{figure}
  
  A look at Figures~\ref{fig:bad} and \ref{fig:switch-cases} shows
  that $\gamma'$ is a quasicycle. Thus, $\sigma$ does not terminate at
  $(i,\infty)$. We need to relate leading hyperedges of $\gamma$ to
  those of $\gamma'$.

  Any leading hyperedge of $\gamma$ that is not incident with $u$ is a
  leading hyperedge of $\gamma'$, and vice versa. We assert that
  neither $au$ nor $bu$ is a leading hyperedge of $\gamma$. If they
  were, they would be redundant (since their size is $2$) and
  $\ptn\pi\infty\infty$ would not be $\pi$-skeletal, contrary to the
  assumption. Finally, if $ucd$ is a leading hyperedge of $\gamma$,
  then $ud$ is a leading hyperedge of $\gamma'$. Note that
  $ucd <_E ud$.

  It follows that if $\sigma$ does not stop at $(i,\infty)$, then
  $d^\pi_i \leq_E d^\sigma_i$. On the other hand, if it does stop,
  then the same inequality holds since $\stopx$ is the greatest
  element of the ordering $\leq_E$. In both cases, we have
  $\pi\worseeqx i\infty\sigma$.

  The last possibility left to consider is $j=0$. We need to show that
  $\ptn\pi i 0 \leq\ptn\sigma i 0$ under the assumption that
  $\pi\eqx {i-1}\infty\sigma$. Let $\RR = \ptn\pi{i-1}\infty$ and let
  $f := d^\pi_{i-1} = d^\sigma_{i-1}$. Since $\ptn\pi\infty\infty$ is
  $\pi$-skeletal, $f\neq\stopx$. If $f = \termx$, then
  $\pi\eq\sigma$. We may thus assume that $f$ is a hyperedge, and in
  that case it is not incident with $u$ (since $H$ and $\tilde H$ have
  no common hyperedge incident with $u$). Thus, $\sigma-f$ is obtained
  from $\pi-f$ by a switch at $u$. Similarly, if $f$ is not used by
  $\pi$, then $\sigma$ is obtained from $\pi$ by a switch at $u$, in the
  hypergraph $H-f$.

  If $f$ is an $X$-antibridge with respect to $\pi$ for some
  $X\in\RR$, we may use Lemma~\ref{l:sub-switch} in the hypergraph
  $H-f$. We find that $\sigma$ is anticonnected on each anticomponent of
  $\pi$ on $X$ in $H-f$, and hence $\pi\worseeqx i 0 \sigma$. On the
  other hand, if $f$ is an $X$-bridge with respect to $\pi$, then
  Lemma~\ref{l:sub-switch} implies that $\sigma-f$ is connected on each
  component of $\pi-f$ on $X$ in $H$, and $\pi\worseeqx i 0 \sigma$
  again. This proves~\eqref{eq:stable-bad} and the lemma follows.
\end{proof}

Let us now state the result we need to use in~\cite{KV-ess-9}.

\begin{theorem}\label{t:no-bad}
  Let $H$ be a $3$-hypergraph. There exists a hypergraph $\tilde H$
  related to $H$ and an acyclic quasigraph $\sigma$ in $\tilde H$ such
  that $\sigma$ has no bad leaves and $V(\tilde H)$ admits a
  $\sigma$-skeletal partition $\SSS$.
\end{theorem}
\begin{proof}
  Let $\MM$ be the set of all quasigraphs in 3-hypergraphs related to
  $H$. Furthermore, let $\MM' \subseteq \MM$ be the set of quasigraphs
  satisfying the following conditions:
  \begin{enumerate}[(1)]
  \item $\sigma$ is $\sqsubseteq$-maximal in $\MM$,
  \item subject to (1), $\sigma$ uses as few hyperedges as possible.
  \end{enumerate}

  Let $\sigma'$ be any quasigraph from $\MM'$; say, $\sigma'$ is a
  quasigraph in a hypergraph $H'$ related to $H$. Theorem
  \ref{t:enhancing} implies that $\sigma'$ is acyclic and
  $\ptn{\sigma'}\infty\infty$ is $\sigma'$-skeletal (where the
  partition is obtained via the plane sequence with respect to
  $H'$). In particular, it makes sense to consider bad leaves of
  $\sigma'$.

  Let us choose $\sigma$ as an element of $\MM'$ with as few bad
  leaves as possible, and define $\SSS := \ptn\sigma\infty\infty$.
  
  We need to prove that $\sigma$ has no bad leaves. Suppose to the
  contrary that there is a bad leaf $u$ for $\sigma$. By
  Lemma~\ref{l:stable-bad}, $\sigma \sqsubseteq \sw\sigma
  u$. Furthermore, $\sw\sigma u$ uses the same number of hyperedges as
  $\sigma$, and has one bad leaf fewer, a contradiction with the
  choice of $\sigma$.
\end{proof}

\section*{Acknowledgment}

We thank two anonymous referees for a careful reading of the
manuscript and many helpful comments.


\end{document}